# On physical diffusion and stochastic diffusion

T. N. Narasimhan

*Although the same mathematical expression is used to describe physical diffusion and stochastic diffusion, there are intrinsic similarities and differences in their nature. A comparative study shows that characteristic terms of physical and stochastic diffusion cannot be placed exactly in one-to-one correspondence. Therefore, judgment needs to be exercised in transferring ideas between physical and stochastic diffusion.*

**Historical context**

Physical diffusion relates to conductive heat transfer in solids, molecular diffusion in substances, flow of galvanic current, and pressure-driven fluid flow in resistive media. It is modelled using the format of the heat equation, introduced by Joseph Fourier in 1807 (ref. 1) in an unpublished monograph. Fourier solved a number of problems relating to solids of symmetrical form and finite dimensions (e.g. rod, prism, cylinder, annulus, cube, sphere). Noteworthy were his solutions involving trigonometric series, specifically designed for bounded domains with the premise that when a function is introduced, its domain is to be specified. At that time, the notion of a function had not yet been rigorously established. Lagrange, one of the four reviewers appointed by the French Academy, erroneously concluded that Fourier's use of trigonometric series had inconsistencies[2], and did not approve the work. The monograph was never published.

Stochastic diffusion relates to spreading associated with random processes. At the turn of the 19th century, Pierre Simon Laplace addressed a fundamental problem of probability theory, propagation of cumulated random values as a function of the number of random values summed up. It was known without formal proof to mathematicians before him that probability distribution of such cumulated values asymptotically approached normal distribution as the number of values summed up became very large. This asymptotic behaviour came to be known as the central limit theorem early during the 20th century at Polyà's suggestion[3]. Between 1809 and 1812, Laplace addressed this problem along two different lines. In 1809, he used génératrice functions[4] to formulate a differential equation of the same mathematical form as Fourier's one-dimensional heat equation, with probability density as the dependent variable, and showed that the equation was satisfied by a normal distribution function. Later, in his masterpiece 'Analytic Theory of Probability', Laplace[5] used characteristic functions, and partially succeeded in proving the central limit theorem[6]. Unlike the bounded domains for which Fourier had designed trigonometric series, Laplace's stochastic diffusion equation pertained to an unbounded, infinite domain, and its solution, the normal distribution function, could not be represented by a trigonometric series.

Influenced by Laplace, Fourier[7] went beyond bounded domains and extended thermal analysis to infinite domains. In doing so, he reproduced Laplace's 1809 result[5] and obtained normal distribution function as a solution to the special case of an instantaneous plane source of heat. By the early 1820s, equations of physical and stochastic diffusion had become established through the contributions of Fourier[8] and Laplace[4]. The equations were unified by normal distribution.

After providing a partial proof of central limit theorem, Laplace did not further pursue stochastic diffusion. Whereas Fourier's physical diffusion became influential, Laplace's stochastic diffusion receded from attention. More than 50 years later, Lord Rayleigh[9,10] and Francis Edgeworth[11] serendipitously discovered connection between normal distribution and Fourier's heat equation, and rejuvenated interest in stochastic diffusion. Luis Bachelier[12] cast, for the first time, the random sampling problem in the form of random variations of stock prices with time. In 1905, Albert Einstein[13] applied stochastic diffusion equation to physically observable random motion of Brownian particles to help establish a molecular-kinetic theory of heat. Immediately following Einstein, Marian Smoluchowski[14] independently arrived at results similar to those of Einstein for Brownian motion. In 1908, Jean Perrin[15] used latest developments in optical microscopy to validate Einstein's theory and confirm molecular nature of matter, a question that had intrigued physics for well over a century. A detailed account of the inter-related history of physical and stochastic diffusion is given in ref. 16.

Soon, two works of major importance added vigour to the investigation of Brownian motion and stochastic diffusion. In 1906, Andrej Markov[17] shifted attention from diffusion equation itself to the process by which the sum $X(n)$ of $n$ random variables changes to $X(n + 1)$ of $(n + 1)$ variables by the addition of a random variable. This work established the study of memory-less Markov processes. Two years later, Paul Langevin[18] provided an alternate mathematical perspective of Brownian motion by expressing the random step of a Brownian particle as a sum of a regular component and a random component. Together, these contributions inspired studies by Norbert Wiener and Andrej Kolmogorov, leading to Kiyoshi Ito's stochastic differential equations.

In the present work, focus is on diffusion equation, examining how it portrays spreading arising from physical and stochastic phenomena. The process governing evolution of a single random event $X(n)$ to $X(n + 1)$ is beyond the scope of this work.

For purposes of comparative study, stochasticity and probability are used here in restricted sense. Conforming to Einstein's assumption that movement of Brownian particles are mutually independent, stochasticity is limited to unrestricted randomness. The second restricted use concerns definition of probability. In studying statistical behaviour of random events, 'relative frequency' is defined as the ratio of number of points falling within a given class to the total number of data points[19]. It is assumed that as the total number of data points becomes large, relative frequency will approximate probability. On mathematical considerations some would posit that relative frequency converges to 'theo-retical probability' as the total number of data points tends to infinity. In this work, attention is restricted to relative frequency, assuming that if the total number of data points is sufficiently large, albeit finite, relative frequency will reasonably represent probability.

The similarities between physical diffusion and stochastic diffusion are fascinating. Fourier (ref. 8, p. 7) himself expressed surprise that the same abstract mathematical expression of pure analysis which occurred in the study of heat diffusion would also arise in all the chief problems of probability theory. Despite mathematical similarities, there are significant differences in the nature of physical and stochastic diffusion. For those in physical and natural sciences, it is of practical interest to examine and appreciate these similarities and differences.

Addressing the relationship between mathematics and physics, James Clerk Maxwell observed[20], 'The figure of speech or of thought by which we transfer the language of a familiar science to one with which we are less acquainted may be called Scientific Metaphor'. Maxwell, who had a flair for poetry, was aware that confusion may arise if metaphor is not appreciated thoughtfully. He went on to add, '. . . The characteristic of a truly scientific system of metaphors is that each term in its metaphorical use retains all the formal relations to the other terms of the system it had in its original use. . . '. In this spirit, this work attempts to identify correspondence between the characteristic terms of the equations describing physical diffusion and stochastic diffusion.

## Temperature, probability density, and conservation

Fourier's one-dimensional heat equation[1] is,

$$K \frac{\partial^2 T}{\partial X^2} = \rho c \frac{\partial T}{\partial t}, \quad (1)$$

where $K$ is thermal conductivity, $T$ the temperature, $\rho$ the mass density, $c$ the specific heat, $X$ the spatial coordinate and $t$ the time. Defining thermal diffusivity, $\eta = K/\rho c$, eq. (1) becomes

$$\eta \frac{\partial^2 T}{\partial X^2} = \frac{\partial T}{\partial t}. \quad (2)$$

Stochastic diffusion equation, originally formulated by Laplace[4], is

$$D \frac{\partial^2 f}{\partial X^2} = \frac{\partial f}{\partial n}, \quad (3)$$

where diffusion coefficient $D$ is analogous to $\eta$, probability density $f(X, n)$ analogous to $T(X, t)$, the sum of $n$ random variables $X$ analogous to spatial coordinate and $n$ analogous to time. Latent in these equations is a conservation principle involving additive extensive quantities.

Temperature $T$ at a point indicates the quantity of heat $\delta H$ stored in an elemental volume $\delta V$ in the point's vicinity and is given by,

$$T = \frac{\delta H}{\rho c \delta V}, \quad (4)$$

where $(\rho c \delta V)$ is heat capacity of $\delta V$. Implicit in eqs (1) and (2) is the fact that heat is conserved over a volume. $T$ is an intensive quantity obtained by normalizing the extensive quantity by a physical property, heat capacity.

In eq. (3), $n$ independent random variables $x$ are assumed drawn from the same distribution with mean $\mu$ and variance $\sigma^2$, and $X(n) = x_1 + x_2 + \ldots + x_n$. Parameter $D$ represents half the variance. Suppose a large number of realizations $N_{\text{Total}} \gg n$ are achieved, and the values of $X(n)$ are organized into a histogram, showing number of realizations $N_{i,n}$ falling in interval between $X_{i-1}$ and $X_i$ where $i = 1, 2, \ldots, I$. We repeat the process by adding randomly drawn $x$ to each $X(n)$. The frequency distribution for $n$ variables and that for $n + 1$ variables will each have its mean and variance. Equation (3) expresses how relative frequency distribution of $X$ evolves from $n$ to $n + 1$ variables. In following this reasoning, density of relative frequency is effectively treated as equivalent to probability density. Note that the evolution of frequency distribution of the same set of realizations is followed as $n$ is progressively increased, and $X$ is allowed to increase without bounds. It is in this sense that stochasticity is limited to unrestricted randomness in this work.

For a given $n$, $f(X + \delta X, n)\delta X = (N_X/N_{\text{Total}})$. Analogous to $T$, $f$ is an intensive quantity. Note that $T$ is not an additive quantity, whereas $f$ is because relative frequency is expressed as a fraction of $N_{\text{Total}}$. This additivity attribute of $f$ has a consequence. Suppose, $N_{i,n}$, $i = 1, 2, 3, \ldots, I$ denotes number of realizations in $i$th interval. Then, $N_{\text{Total}} = \Sigma_{i=1}^{i=I} N_{i,n}$, and effective probability associated with $i$th interval is equal to $N_{i,n}/N_{\text{Total}}$.

Now suppose $N$ additional realizations are added in the interval $i = j$. This will increase $N_{\text{Total}}$ to $(N_{\text{Total}} + N)$, and change effective probability of interval $j$ to

$$\frac{N_{j,n} + N}{N_{\text{Total}} + N}.$$

Additionally, effective probability of all other intervals $i \neq j$ will change to

$$\frac{N_{i,n}}{N_{\text{Total}} + N}.$$

This implies that the act of changing $N_{\text{Total}}$ induces a change in effective probability throughout the frequency distribution. This change has occurred without the addition of any random variable. However, stochastic diffusion eq. (3) accounts only for changes in probability density arising from addition of random variables, but not change induced by increasing the number of realizations.

Now consider conservation. In the heat problem, quantity of heat $H$, an extensive quantity, is conserved over a domain. Analogous to this, in the stochastic problem, number of realizations is an extensive quantity and is to be conserved. In the case of an instantaneous heat source, the quantity of heat released at $X = 0$ at $t = 0$ is analogous to $N_{\text{Total}}$. This analogy implies that stochastic diffusion eq. (3) expresses the evolution of a system driven purely by initial conditions over an unbounded domain, solely due to random changes.

## Dimensionality

As introduced by Laplace, stochastic differential equation has only one independent variable associated with the second derivative. Therefore, eq. (3) is analogous to Fourier's one-dimensional heat equation (eq. (2)). An interesting question concerns whether eq. (3) can be extended to more than one dimension. Rayleigh[10] addressed this issue in his study of mixing random acoustic waves of the same amplitude, but of arbitrary phase. First he considered a problem in which phases were randomly chosen to be +1 or –1. This distribution has vari-ance $\sigma^2 = 1$. For this case, he obtained the differential equation,

$$\frac{1}{2}\frac{\partial^2 f}{\partial X^2} = \frac{\partial f}{\partial n}. \qquad (5)$$

He then extended analysis to phases distributed over the interval from –1 to +1. In this case, the resultant of $n$ random vibrations had an amplitude $r$ and phase $\theta$. Representing $(r, \theta)$ in Cartesian coordinates, he derived two-dimensional form,

$$\frac{1}{4}\left[\frac{\partial^2 f}{\partial X^2} + \frac{\partial^2 f}{\partial Y^2}\right] = \frac{\partial f}{\partial n}. \qquad (6)$$

Going beyond theory of vibrations, Rayleigh[10] found that the method could be extended to composition of unit vectors in three dimensions whose directions were randomly chosen. For this case, the three-dimensional stochastic equation became

$$\frac{1}{6}\left[\frac{\partial^2 f}{\partial X^2} + \frac{\partial^2 f}{\partial Y^2} + \frac{\partial^2 f}{\partial Z^2}\right] = \frac{\partial f}{\partial n}. \qquad (7)$$

With the application of stochastic diffusion to Brownian motion, dimensionality is determined by physical dimensions and time.

## Conductivity, capacity and diffusion coefficient

As seen from eqs (2) and (3), thermal diffusivity is analogous to one-half of variance. Thermal diffusivity, however, is a function of two fundamentally important parameters, thermal conductivity and heat capacity, each amenable to accurate measurement with instruments. In resistive materials that host compressible fluid flow, conductivity is related to the material's hydraulic resistance, and hydraulic capacitance (analogous to heat capacity) is related to the fluid's compressibility. Transient diffusion of heat or fluid is governed by an interaction between conductivity, capacitance and spatial variation of appropriate intensive quantity, temperature or fluid pressure. In stochastic diffusion, separate attributes analogous to conductivity and capacitance are not identifiable.

## Bounded domains and stochasticity

Heat diffusion is driven by spatial variations in temperature. Such spatial variations may be those that exist at time $t = 0$ (initial conditions), or may arise due to conditions imposed on the boundaries of finite domains. In infinite domains, diffusion can arise due to initial conditions, or due to sources or sinks. Sources (plane, line or point) may be instantaneous or continuous. We have already seen that because of the way relative frequency is defined, stochastic diffusion equation (eq. (3)) is restricted to problems with constant $N_{Total}$, that is, problems analogous to instantaneous sources, driven purely by initial conditions.

It may be argued that this constraint can be overcome by using theoretical probability with $N_{Total}$ tending to infinity. If so, probability density, an intensive quantity, will also serve as the extensive quantity that is conserved. A consequence will be that stochastic diffusion and physical diffusion will notably differ in regard to conservation. In physical diffusion, intensive and extensive quantities are mutually distinct, and are related through well-defined empirical relations. Consequently, in solving a physical diffusion equation, it is important to assure that the function sought as a solution also satisfies a separate conservation criterion. In numerical simulations, conservation helps assure credibility of solutions. If the dependent variable, an intensive quantity, also serves as the extensive quantity that is conserved, the situation would significantly differ from that of physical diffusion. The role of conservation in stochastic diffusion merits further study.

In bounded domains, heat diffusion is influenced by conditions (temperature or heat flux) prescribed on the boundary. In Laplace's stochastic problem, $X(n)$, imitates spatial dimension. Therefore, a 'bounded' domain implies that $X(n)$ for any $n$ will have to be within prescribed values. Suppose the domain is bounded by $-X_A < X(n) > X_A$. If, the addition of a random variable causes $X(n + 1)$ to fall below $-X_A$ or exceed $X_A$, two courses of action may be possible. The realization may be rejected because it falls outside the bound, and a new realization added to keep $N_{Total}$ constant. Or, the value of $X(n + 1)$ maybe modified in some way to keep the value within bounds, and keep $N_{Total}$ constant. Either action would introduce bias and inhibit randomness.

Spreading in an infinite domain caused by initial conditions, and equilibrium configuration of boundary-value problems have distinct mathematical attributes. Normal distribution of the former cannot be represented by trigonometric series. In fact, Lagrange's criticism of Fourier's work was that trigonometric series could not be extended to unbounded domains[2]. In comparison, steady-state configurations of boundary-value problems are governed by variational calculus and Dirichlet principle. By definition, solutions to boundary-value problems are harmonic functions.

## Brownian motion

So far, stochastic diffusion has been viewed in terms of abstract thought experiments. In Brownian motion, randomness is associated with an observable physical process, and stochastic diffusion can be physically understood in three dimensions. Here, the number of realizations $N_{Total}$ is supplanted by the total number of Brownian particles released at $x = 0$ at time zero. Thus $f(X, t)dX$ denotes $(N_{X,t}/N_{Total})$, where $N_{X,t}$ is the number of particles in the immediate vicinity of $X$. Conservation implies that one follows the position of the same set of particles as time progresses.

Einstein[13] presented the following equation for colloidal particles in a state thermally agitated perpetual motion in water.

$$D_{st}\frac{\partial^2 f}{\partial X^2} = \frac{\partial f}{\partial t}, \qquad (8)$$

where $f(X, t)dX$ is the probability that a particle lies in the immediate vicinity of $X$ at time $t > 0$. The stochastic diffusion coefficient $D_{st}$ is given by,

$$D_{st} = \frac{\sigma^2}{2\tau}, \qquad (9)$$

where $\tau$ is time interval, sufficiently small in comparison with $t$, but sufficiently large so that movements executed by a particle over two successive time intervals are mutually independent. For particles of diameter 1 micron in water at 17°C, $\tau$ is of the order of $10^{-8}$ s (ref. 13).

Note that if $N_{Total}$ becomes very large, Brownian particles will collide among themselves, cause flow resistance to increase, and affect rate of diffusion. For this reason, Einstein restricted attention to dilute suspensions and independent random movements. Therefore, total number of particles, $N_{Total}$, cannot be allowed to tend to infinity to define a theoretical probability.

Contemporaneously with Einstein, Karl Pearson[21] introduced the concept of random walk in which a man starts from a point $O$ and walks $\ell$ yards in a straight line, then turns through any angle whatever and walks another $\ell$ yards. This process is repeated $n$ times. The random-walk problem is to estimate the probability that after $n$ stretches, the man is at a distance between $r$ and $r + \delta r$ from the starting point $O$. Viewed in this context, the Brownian motion problem consists in a particle starting from $O$ moving in some direction for an interval of time $\delta t = \tau$, then turns through any angle whatever and moving for another interval of time $\tau$, and repeating the process $n$ times. This implies that total time $t = n\tau$, and $dt = \tau\,dn$. Consequently, Einstein's equation (eq. (8)) reduces to

$$\tau D_{st} \frac{\partial^2 f}{\partial X^2} = \frac{\partial f}{\partial n}. \quad (10)$$

However, in view of eq. (10), $\tau D_{st} = \sigma^2/2 = D$. Therefore eq. (10) is essentially the same as Laplace's equation for random variables.

With admirable insight, Einstein saw in Brownian motion both a micro-scale stochastic process and a macro-scale non-random deterministic process. On the microscopic scale, the motion of colloidal particles was an unrestricted random process. Macroscopically, it was a non-random deterministic process of molecular diffusion amenable to precise description using Fick's Law[22]. Using a dynamical form of Fick's Law by balancing impelling osmotic forces and resisting hydraulic forces, Einstein obtained a macroscopic diffusion coefficient $D_M$ for spherical colloidal particles diffusing in viscous water. Osmotic pressure, in turn, was a colligative property, solely depend on number of suspended particles[23]. Thus,

$$D_M = \frac{RT}{6 N_{Avo} \pi \mu r}, \quad (11)$$

where, $R$ is universal gas constant, $T$ the absolute temperature, $\mu$ the coefficient of viscosity, $N_{Avo}$ the Avogadro's number and $r$ the radius of a spherical particle. In suggesting an experimental procedure for establishing molecular-kinetic theory of heat, Einstein's central assumption was that $D_{st}$ (an attribute of a frequency distribution), and $D_M$ (a physically measurable quantity) were equal in magnitude and that they could be used interchangeably.

## Concluding remark

Physical and stochastic diffusion are widely used in physical, biological and earth sciences and engineering. The present work compares and contrasts elements of physical diffusion on the one hand and stochastic diffusion on the other so that ideas can be credibly exchanged between the two fields. The study shows that not all elements of the equations of physical and stochastic diffusion can be placed exactly in one-to-one correspondence. Consequently, thought and judgement are necessary in transferring ideas among these models.

ACKNOWLEDGEMENTS. I am thankful to James Pitman for illuminating discussions on the basics of probability theory. Thanks are due to Andre Journel for reviewing the draft of the manuscript. This work was partly supported by an enabling grant from the Committee on Research, University of California at Berkeley.

*T. N. Narasimhan is in the Department of Materials Science and Engineering, 210 Hearst Memorial Mining Building, University of California, Berkeley, CA 94720-1760, USA.*
*e-mail: tnnarasimhan@LBL.gov*